\documentclass[fleqn]{article}

\usepackage{amsmath}
\usepackage{amsfonts}
\usepackage{xcolor}
\usepackage{calrsfs}

\newcommand{\FF}{\mathbb{F}}

\newcommand{\PG}{\mathrm{PG}}

\newcommand{\cQ}{{\cal Q}}
\newcommand{\cH}{{\cal H}}
\newcommand{\cS}{{\cal S}}
\newcommand{\cM}{{\cal M}}
\newcommand{\ve}{\varepsilon}

\newcommand{\veS}{\varepsilon_{\mathrm{nat}}}
\newcommand{\Pu}{\mathrm{Pure}}
\newcommand{\ovF}{\overline{\FF}}
\newcommand{\ovV}{\overline{V}}

\newcommand{\bx}{{\bf x}}
\newcommand{\by}{{\bf y}}

\newcommand{\be}{{\bf e}}
\newcommand{\ba}{{\bf a}}
\newcommand{\bfb}{{\bf b}}

\newcommand{\bv}{{\bf v}}

\newcommand{\gA}{\mathfrak{a}}
\newcommand{\gB}{\mathfrak{b}}
\newcommand{\gC}{\mathfrak{c}}
\newcommand{\gX}{\mathfrak{x}}

\newtheorem{conc}{Conclusion}[section]

\newtheorem{theo}[conc]{Theorem}
\newtheorem{cor}[conc]{Corollary}
\newtheorem{prop}[conc]{Proposition}
\newtheorem{lemma}[conc]{Lemma}

\newtheorem{note}{Remark}

\title{Geometric hyperplanes of the Lie geometry $A_{n,\{1,n\}}(\FF)$} 
\author{Antonio Pasini}
\date{} 

\begin{document}

\maketitle

\begin{abstract} 
In this paper we investigate hyperplanes of the point-line geometry $A_{n,\{1,n\}}(\FF)$ of point-hyerplane flags of the projective geometry $\PG(n,\FF)$. Renouncing a complete classification, which is not yet within our reach, we describe the hyperplanes which arise from the natural embedding of $A_{n,\{1,n\}}(\FF)$, that is the embedding which yields the adjoint representation of $\mathrm{SL}(n+1,\FF)$. By exploiting properties of a particular sub-class of these hyerplaes, namely the {\em singular hyperplanes}, we shall prove that all hyperplanes of $A_{n,\{1,n\}}(\FF)$ are maximal subspaces of $A_{n,\{1,n\}}(\FF)$. Hyperplanes of $A_{n,\{1,n\}}(\FF)$ can also be contructed starting from suitable line-spreads of $\PG(n,\FF)$ (provided that $\PG(n,\FF)$ admits line-spreads, of course). Explicitly, let $\mathfrak{S}$ be a composition line-spread of $\PG(n,\FF)$ such that every hyperplane of $\PG(n,\FF)$ contains a sub-hyperplane of $\PG(n,\FF)$ spanned by lines of $\mathfrak{S}$. Then the set of points $(p,H)$ of $A_{n,\{1,n\}}(\FF)$ such that $H$ contains the member of $\mathfrak{S}$ through $p$ is a hyperplane of $A_{n,\{1,n\}}(\FF)$. We call these hyperplanes {\em hyperplanes of spread type}. Many but not all of them arise from the natural embedding. 
\end{abstract} 

\section{A selection of results}

\subsection{Basic properties of the geometry $A_{n,\{1,n\}}(\FF)$}

Following a well established notation, we denote by $A_{n,\{1,n\}}(\FF)$ the geometry of point-hyperplane flags of the projective geometry $\PG(n,\FF)$, for $n \geq 2$ and $\FF$ a given field, namely a commutative division ring. 

Explicitly, $A_{n,\{1,n\}}(\FF)$ is the point-line geometry the points of which are the ordered pairs $(p,H)$ where $p$ and $H$ are a point and a hyperplane of $\PG(n,\FF)$ respectively and $p\in H$; the lines of $A_{n,\{1,n\}}(\FF)$ are the sets $\{(p,H)\mid p\in \ell\}$ for $\ell$ a line of $\PG(n,\FF)$ and $H$ a hyperplane of $\PG(n,\FF)$ containing $\ell$ and the sets $\{(p,H)\mid H \supset L\}$ for $L$ a sub-hyperplane of $\PG(n,\FF)$ (namely a subspace of $\PG(n,\FF)$ of codimension 2) and $p$ a point of $L$. Accordingly, two points $(p,H)$ and $(q,K)$ of $A_{n,\{1,n\}}(\FF)$ are collinear if and only if either $p = q$ or $H = K$. If $p\neq q$ and $H\neq K$ but either $p\in K$ or $q\in H$, then $(p,H)$ and $(q,K)$ are at distance 2, otherwise they are at distance 3. Thus, the diameter of the collinearity graph of $A_{n,\{1,n\}}(\FF)$ is equal to 3. 

\subsubsection{Polar and special pairs}

When $n = 2$ the point-line geometry $A_{n,\{1,n\}}$ is a generalized hexagon where every point belongs two just two lines. 

Let $n > 2$. Then $A_{n,\{1,n\}}(\FF)$ is a parapolar space of symplectic rank 2, as defined by Shult \cite{S}, the symps of which are grids isomorphic to the hyperbolic quadric $\cQ^+_3(\FF)$. Recall that a pair of points $\{X,Y\}$ at distance 2 of a parapolar space is called {\em polar} if it is contained in a symp, equivalently at least two points exist which are collinear with both $X$ and $Y$. Otherwise, when just one point exists which is collinear with both $X$ and $Y$, the pair $\{X,Y\}$ is said to be {\em special}. A pair $\{(p,H), (q,K)\}$ of points of $A_{n,\{1,n\}}(\FF)$ at distance 2 is polar precisely when $p\in K$ and $q\in H$. If this is the case, the points of the unique symp containg $(p,H)$ and $(q,K)$ are the point-hyperplane flags $(x,X)$ of $\PG(n,\FF)$ with $x$ a point of the line $\ell = \langle p,q\rangle$ and $X$ a hyperplane of $\PG(n,\FF)$ containing the sub-hyperplane $H\cap K$. On the other hand, let $\{(p,H),(q,K)\}$ be a special pair, say $p\in K$ but $q\not\in H$, to fix ideas. Then $(p,K)$ is the unique point of $A_{n,\{1,n\}}(\FF)$ collinear with both $(p,H)$ and $(q,K)$. 

\subsubsection{Maximal singular subspaces} 

Given a point $a$ of $\PG(n,\FF)$, let ${\cal M}_a$ be the set of pairs $(a,H)$ with $H$ a hyperplane of $\PG(n,\FF)$ containing $a$. This set is a maximal singular subspace of $A_{n,\{1,n\}}(\FF)$. Dually, for a hyperplane $A$ of $\PG(n,\FF)$, the set ${\cal M}_A := \{(p,A)\mid p\in A\}$ is a maximal singular subspace of $A_{n,\{1,n\}}(\FF)$. We say that ${\cal M}_a$ and ${\cal M}_A$ are {\em based} at $a$ and $A$ respectively. Every maximal singular subspace of $A_{n,\{1,n\}}(\FF)$ admits one of these two descriptions. So, the maximal singular subspaces of $A_{n,\{1,n\}}(\FF)$ are partioned in two families: those which are based at a point of $\PG(n,\FF)$ and those based at a hyperplane of $\PG(n,\FF)$. Two distinct maximal singular subspaces have at most one point in common; they meet in a point only if they do not belong to the same family. Moreover. if $M$ is a maximal singular subspace of $A_{n,\{1,n\}}(\FF)$ and $(p,H)$ a point of $A_{n,\{1,n\}}(\FF)$ exterior to $\cal M$, then $(p,H)$ is collinear with at most one point of $\cal M$. 

\subsubsection{The natural embedding}          

It is well known that $A_{n,\{1,n\}}(\FF)$ admits a (full) projective embedding in the projective space $\PG(M_{n+1}^0(\FF))$ of the vector space $M_{n+1}^0(\FF)$ of null-traced square matrices of order $n+1$ with entries in $\FF$, which yields the adjoint representation of the special linear group $\mathrm{SL}(n+1,\FF)$. We call it the {\em natural embedding} of $A_{n,\{1,n\}}(\FF)$ and we denote it by the symbol $\veS$. 

Explicitly, recall that $M^0_{n+1}(\FF)$ is a hyperplane of the vector space $M_{n+1}(\FF)$ of square matrices of order $n+1$ with entries in $\FF$ and the latter is the same as the tensor product $V\otimes V^*$, where $V = V(n+1,\FF)$ and $V^*$ is the dual of $V$. The pure tensors $\bx\otimes \xi$ of $V\otimes V^*$, with $\bx$ and $\xi$ non-zero vectors of $V$ and $V^*$ respectively, yield the matrices of $M_{n+1}(\FF)$ of rank $1$.  With $\bx$ and $\xi$ as above, let $[\bx]$ and $[\xi]$ be the point and the hyperplane of $\PG(n,\FF)$ represented by $\bx$ and $\xi$. Then $([\bx], [\xi])$ is a point of $A_{n,\{1,n\}}(\FF)$ if and only if $\xi(\bx) = 0$. The pure tensor $\bx\otimes \xi$, regarded as a square matrix of $M_{n+1}(\FF)$ of rank 1, is null-traced if and only if $\xi(\bx) = 0$. The natural embedding 
\[\veS :A_{n,\{1,n\}}(\FF)\rightarrow\PG(M^0_{n+1}(\FF))\]
maps the point $([\bx],[\xi])$ of $A_{n,\{1,n\}}(\FF)$ onto the point $\langle \bx\otimes \xi\rangle$ of $\PG(M^0_{n+1}(\FF))$. These points form a projective variety, commonly known as the Segre variety. 

\subsection{Generalities on hyperplanes of point-line geometries}  

In this subsection we recall some definitions and basics facts on subspaces and geometric hyperplanes of point-lines geometries, to be exploited in the sequel. 

Throughout this subsection $\Gamma$ is an arbitrary point-line geometry as defined in Shult \cite{S}, but we assume that the lines of $\Gamma$ are subsets of the point-set of $\Gamma$. A {\em subspace} of $\Gamma$ is a subset $X$ of the point-set of $\Gamma$ such that, for every line $\ell$ of $\Gamma$, if $|\ell\cap X| > 1$ then $\ell\subseteq X$. A proper subspace of $\Gamma$ is said to be a {\em geometric hyperplane} of $\Gamma$ (a {\em hyperplane} of $\Gamma$ for short) if every line of $\Gamma$ meets it non-trivially. Equivalently, a hyperplane of $\Gamma$ is a proper subset $\cH$ of the point-set of $\Gamma$ such that, for every line $\ell$ of $\Gamma$, either $\ell \subseteq \cH$ or $|\ell\cap \cH| = 1$. 

Assuming that $\Gamma$ is embeddable, let $\ve:\Gamma\rightarrow\Sigma$ be an embedding of $\Gamma$ in a projective space $\Sigma$, as defined by Shult in \cite{S93}. We say that a hyperplane $\cH$ of $\Gamma$ {\em arises} from $\ve$ when $\ve(\cH)$ spans a projective hyperplane $H$ of $\Sigma$ and $\cH = \ve^{-1}(H)$. Conversely, the $\ve$-preimage $\ve^{-1}(H)$ of a projective hyperplane $H$ of $\Sigma$ is a hyperplane of $\Gamma$ and $\ve(\ve^{-1}(H)) = H\cap\ve(\Gamma)$, where $\ve(\Gamma)$ is the $\ve$-image of the point-set of $\Gamma$. Hwever in general $\ve^{-1}(H)$ does not arise from $\ve$. Indeed $\ve^{-1}(H)$ arises from $\ve$ if and only if $\langle H\cap\ve(\Gamma)\rangle = H$. 

\begin{prop}\label{Gen 1} 
Given a projective hyperplane $H$ of $\Sigma$, suppose that $\ve^{-1}(H)$ is maximal in the family of proper subspaces of $\Gamma$. Then $\langle H\cap\ve(\Gamma)\rangle = H$. 
\end{prop}
{\bf Proof.} For a contradiction, let $\langle H\cap\ve(\Gamma)\rangle \subset H$. Let $p$ be a point of $\Gamma$ not in $\cH := \ve^{-1}(H)$. Then $\ve(p)\not\in H$ and, since $\ve(\cH) = H\cap\ve(\Gamma)$ does not span $H$, the set $\ve(\cH)\cup\{\ve(p)\}$ does not span $\Sigma$. However, $\cH\cup\{p\}$ generates $\Gamma$, since $\cH$ is a maximal subspace of $\Gamma$ by assumption, and $\ve(\Gamma)$ spans $\Sigma$ (by the definition of embedding). Hence $\ve(\cH)\cup\{\ve(p)\}$ spans $\Sigma$. Contradiction.  \hfill $\Box$  

\subsection{The hyperplanes of $A_{n,\{1,n\}}(\FF)$ which arise from $\veS$}\label{tensor}  

Let $f : M_{n+1}(\FF)\times M_{n+1}(\FF) \rightarrow \FF$ be the symmetric bilinear form defined as follows: for any two matrices $X,Y \in M_{n+1}(\FF)$, we put
\begin{equation}\label{f1}
f(X, Y) ~ := ~ \mathrm{trace}(XY)
\end{equation}
where $XY$ is the usual row-times-column product. More explicitly, with $X = (x_{i,j})_{i,j=0}^n$ and $Y = (y_{i,j})_{i,j=0}^n$, we have  
\begin{equation}\label{f2} 
f((x_{i,j})_{i,j=0}^n, (y_{i,j})_{i,j=0}^n) ~ = ~ \sum_{i,j}x_{i,j}y_{j,i}.
\end{equation}
In particular, when $X$ and $Y$ have rank 1, namely they are pure tensors, say $X = \bx\otimes\xi$ and $Y = \by\otimes \upsilon$, then
\begin{equation}\label{f3} 
f(\bx\otimes\xi, \by\otimes\upsilon) ~ = ~ \upsilon(\bx)\xi(\by).
\end{equation}
In tensor language, with a matrix $X = (x_{i,j})_{i,j=0}^n$ regarded as a tensor $x^i_{.j}$, the scalar $\mathrm{trace}(XY)$ is the complete saturation of the tensors $X$ and $Y$. So, we call $f$ the {\em saturation form} of $M_{n+1}(\FF)$. 

The saturation form $f$ is non-degenerate (and trace-valued when $\mathrm{char}(\FF) \neq 2$). Therefore, if $\perp_f$ is the orthogonality relation associated to $f$, the hyperplanes of $M_{n+1}(\FF)$ are the sets $M^{\perp_f}$ for $M\in M_{n+1}(\FF)\setminus\{O\}$ (where $O$ stands for the null matrix). In particular, $M^{\perp_f} = M_{n+1}^0(\FF)$ if and only if $M$ is proportional to the identity matrix $I$ (and $M\neq O$, of course). For $M\in M_{n+1}(\FF)\setminus\{\lambda I\}_{\lambda\in \FF}$, we put
\[\cH_M ~ := ~ \veS^{-1}([M^{\perp_f}\cap M_{n+1}^0(\FF)])\]
where $[M^{\perp_f}\cap M_{n+1}^0(\FF)]$ is the hyperplane of $\PG(M^0_{n+1}(\FF))$ corresponding to the hyperplane $M^{\perp_f}\cap M_{n+1}^0(\FF)$ of $M_{n+1}^0(\FF)$. So, $\cH_M$ is a geometric hyperplane of $A_{n,\{1,n\}}(\FF)$. We call it a {\em hyperplane of plain type}. 

As we shall see later (Corollary \ref{main2 bis}), the hyperplanes of plain type are precisely the hyperplanes of $A_{1,\{1,n\}}(\FF)$ whch arise from $\veS$.  

\begin{note}\label{rem notation}
\em
Our notation $\cH_M$ conceals the fact that the mapping $M\mapsto \cH_M$ is not injective. Explicitly, for $M, N\in  M_{n+1}(\FF)\setminus\{\lambda I\}_{\lambda\in \FF}$, we have $\cH_M = \cH_N$ if and only if the pairs $\{M,I\}$ and $\{N,I\}$ span the same subspace of $M_{n+1}(\FF)$. The hyperplanes $\cH_M$ of $A_{n,\{1,n\}}(\FF)$ are indeed the points of a projectve geometry isomorphic to the star of the point $\langle I\rangle = \{\lambda I\}_{\lambda\in \FF}$ in $\PG(M_{n+1}(\FF))$. Of course, nothing forbids us to modify our notation in such a way as to keep a record of these facts in it but the resulting notation would be not so nice. We prefer to stick to our slightly incorrect notation. 
\end{note} 
   
The next proposition will be proved in Section \ref{pure}

\begin{prop}\label{hyp pure} 
The following are equivalent for a hyperplane $\cH_M$ of $A_{n,\{1,n\}}(\FF)$ of plain type: 
\begin{itemize}
\item[$(1)$] the hyperplane $\cH_M$ contains a maximal singular subspace of $A_{n,\{1,n\}}(\FF)$ based at point of $\PG(n,\FF)$;
\item[$(2)$] the hyperlane  $\cH_M$ contains a maximal singular subspace of $A_{n,\{1,n\}}(\FF)$ based at a hyperplane of $\PG(n,\FF)$;
\item[$(3)$] the matrix $M$ admits at least one eigenvalue in $\FF$. 
\end{itemize}
\end{prop}

In fact we shall prove the following sharper statement. For a point $a$ of $\PG(n,\FF)$ let $\ba\in V$ a representative vector of $a$. Then $\cH_M$ contains the maximal singular subspace ${\cal M}_a$ if and only if $\ba$ is a (right) eigenvector of $M$, namely $M\ba\in \langle \ba\rangle$. Dually, for a hyperplane $A$ of $\PG(n,\FF)$, let $\alpha\in V^*$ represent $A$. Then $\cH_M\supseteq {\cal M}_A$ if and only if $\alpha M \in \langle\alpha\rangle$ (namely $\alpha$ is a left eigenvector of $M$). Since $M$ admits right eigenvectors (in $V$) if and only if it admits left eigenvectors (in $V^*)$ if and only if it admits eigenvalues in $\FF$, Propositon \ref{hyp pure} follows.   

In particular, let $\mathrm{rank}(M) = 1$, say $M = \ba\otimes \alpha$ for a nonzero vector $\ba\in V$ and a non-zero linear functional $\alpha \in V^*$. Then $M\ba = \ba\cdot \alpha(\ba)$ and $\alpha M = \alpha(\ba)\cdot \alpha$, namely $\ba$ and $\alpha$ are respectively a right eigenvector and a left eigenvector of $M$ and $\alpha(\ba)$ is the eigenvalue associated with both $\ba$ anf $\alpha$. Accordingly, $\cH_M = \cH_{\ba\otimes \alpha}$ contains both ${\cal M}_a$ and ${\cal M}_A$, where $a$ and $A$ are the point and the hyperplane of $\PG(n,\FF)$ represented by $\ba$ and $\alpha$ respectively. 

As noticed in Remark \ref{rem notation}, if $\ba\otimes\alpha$ and $\bfb\otimes\beta$ define the same hyperplane then $\bfb\otimes\beta \in\langle \ba\otimes\alpha, I\rangle$. However, all rank 1 matrices contained in $\langle \ba\otimes\alpha, I\rangle$ are proportional to $\ba\otimes\alpha$. Hence $\cH_{\ba\otimes\alpha} = \cH_{\bfb\otimes\beta}$ if and only if $\ba\otimes\alpha$ and $\bfb\otimes\beta$ are proportional, namely they correspond to the same pair $(a, A)$, where $a$ and $A$ are the point and the hyperplane of $\PG(n,\FF)$ represented by $\ba$ and $\alpha$ respectively. In view of this fact, henceforth we write $\cH_{a,A}$ instead of $\cH_{\ba\otimes\alpha}$. The following will be proved in Section \ref{pure}.

\begin{prop}\label{main1}
Let $a$ and $A$ be a point and a hyperplane of $\PG(n,\FF)$. Then $\cH_{a,A}$ is the set of points of $A_{n,\{1,n\}}(\FF)$ which are collinear with at least one point of $\cM_a\cup\cM_A$. In particular, if $a\in A$ then $(a,A)$ is a point of $A_{n,\{1,n\}}(\FF)$ and $\cH_{a,A}$ is the set of points of $A_{n,\{1,n\}}(\FF)$ at non-maximal distance from $(a,A)$.   
\end{prop}

We call the hyperplanes as $\cH_{a,A}$ {\em quasi-singular} hyperplanes. When $a\in A$ we say that the hyperplane $\cH_{a,A}$ is {\em singular} with $(a,A)$ as its {\em deepest point}.

\subsection{All hyperplanes of $A_{n,\{1,n\}}(\FF)$ are maximal subspaces}

The following will be proved in Section \ref{max}.

\begin{lemma}\label{main2 lem}
All singular hyperplanes of $A_{n,\{1,n\}}(\FF)$ are maximal subspaces of $A_{n,\{1,n\}}(\FF)$.
\end{lemma}

In Section \ref{max}, with the help of Lemma \ref{main2 lem}, we shall prove the following.  

\begin{theo}\label{main2}
All hyperplanes of $A_{n,\{1,n\}}(\FF)$ are maximal subspaces of the point-line geometry $A_{n,\{1,n\}}(\FF)$.
\end{theo} 

Therefore, by Shult \cite[Lemma 4.1.1]{S},

\begin{cor}\label{main2 cor}
Let $\cH$ be a hyperplane of $A_{n,\{1,n\}}(\FF)$. Then the collinearity graph of $A_{n,\{1,n\}}(\FF)$ induces a connected graph on the complement of $\cH$. 
\end{cor}

By combining Theorem \ref{main2} with Proposition \ref{Gen 1}, we immediately obtain the following:

\begin{cor}\label{main2 bis} 
The hyperplanes of $A_{n,\{1,n\}}(\FF)$ which arise from the natural embedding $\veS$ are precisely those of plain type. 
\end{cor}    

\subsection{Hyperplanes of $A_{n,\{1,n\}}(\FF)$ of spread type}\label{intro spread}  

Recall that a {\em line-spread} of $\PG(n,\FF)$ is a family $\mathfrak{S}$ of lines of $\PG(n,\FF)$ such that every point of $\PG(n,\FF)$ belongs to exactly one member of $\mathfrak{S}$. A {\em dual line-spread} of $\PG(n,\FF)$ is a line-spread of the dual of $\PG(n,\FF)$. The lines of the dual of $\PG(n,\FF)$ are in fact sub-hyperplanes of $\PG(n,\FF)$. So, a dual line-spread of $\PG(n,\FF)$ is a family $\mathfrak{S}^*$ of sub-hyperplanes of $\PG(n,\FF)$ such that every hyperplane of $\PG(n,\FF)$ contains exactly one member of $\mathfrak{S}^*$. (Note that $n = 1$ is allowed in these definitions: a $1$-dimensional projective space is a line, say $\ell$, and we regard $\{\ell\}$ as the unique line-spread of $\ell$.) 

As $\PG(n, \FF)$ and its dual are isomorphic (since $\FF$ is a field), if $\PG(n,\FF)$ admits line-spreads then it also admits dual line-spreads. We recall that $\PG(n,\FF)$ admits a line-spread only if $n$ is odd. 
   
When $\PG(n,\FF)$ admits line-spreads, an intriguing family of hyperplanes of $A_{n,\{1,n\}}(\FF)$ also exists, to be described in a few lines (see below, Theorem \ref{main2}). We firstly state a couple of lemmas, to be proved in Section \ref{prel-spread-4}.   

\begin{lemma}\label{main3 lem1} 
Let $\mathfrak{S}$ and $\mathfrak{S}^*$ be a line-spread and a dual line-spread of $\PG(n,\FF)$. The following two properties are equivalent 
\begin{itemize}
\item[$(\mathrm{S})$] for every $L\in\mathfrak{S}^*$, the members of $\mathfrak{S}$ contained in the sub-hyperplane $L$ form a line-spread of the projective space $L$.
\item[$(\mathrm{S}^*)$] for every $\ell\in\mathfrak{S}$, the members of $\mathfrak{S}^*$ containing the line $\ell$ form a dual line-spread in the star of $\ell$.  
\end{itemize}
\end{lemma} 

\begin{lemma}\label{main3 lem2}
For every line-spread $\mathfrak{S}$ of $\PG(n,\FF)$, there is at most one dual line-spread $\mathfrak{S}^*$ such that the pair $\{\mathfrak{S}, \mathfrak{S}^*\}$ satisfies property $(\mathrm{S})$ (equivalently, $(\mathrm{S}^*)$) of Lemma \ref{main3 lem1}.   
\end{lemma} 

Of course, the dual of Lemma \ref{main3 lem2} also holds true: for every dual line-spread $\mathfrak{S}^*$, there exists at most one line-spread $\mathfrak{S}$ such that $\{\mathfrak{S}, \mathfrak{S}^*\}$ satifies $(\mathrm{S})$. 

For a line-spread $\mathfrak{S}$ of $\PG(n,\FF)$, if a dual line-spread $\mathfrak{S}^*$ exists such that $\{\mathfrak{S}, \mathfrak{S}^*\}$ satisfies  property $(\mathrm{S})$ of Lemma \ref{main3 lem1}, we say that $\mathfrak{S}$ {\em admits the dual} and we call $\mathfrak{S}^*$ (unique by Lemma \ref{main3 lem2}) the {\em dual} of $\mathfrak{S}$. 

\begin{note}\label{composition}
\em Recall that a {\em composition line-spread} is a line-spread $\mathfrak{S}$ such that, for every nonempty subset $X$ of $\mathfrak{S}$, the lines of $\mathfrak{S}$ contained in the span $\langle \cup_{\ell\in X}\ell\rangle$ of the join of the members of $X$ form a line-spread of that span. It can be proved (Van Maldeghem \cite{Hpc1}) that a line-spread $\mathfrak{S}$ of $\PG(n,\FF)$ admits the dual if and only if $\mathfrak{S}$ is a composition line-spread and every hyperplane of $\PG(n,\FF)$ contains a sub-hyperplane of $\PG(n,\FF)$ spanned by lines of $\mathfrak{S}$, the family $\mathfrak{S}^*$ of these sub-hyperplanes being the dual of $\mathfrak{S}$. In particular, when $\FF$ is finite $\mathfrak{S}$ admits the dual if and only if it is a composition line-spread.  However we are not going to prove these facts in this paper. Property $(\mathbb{S})$ is sufficient for our purposes in this paper.   
\end{note}   

The next lemma will also be proved in Section \ref{prel-spread-4}. 

\begin{lemma}\label{main3 lem3}
Let $\mathfrak{S}$ be a line-spread of $\PG(n,\FF)$ admitting the dual and let $\mathfrak{S}^*$ be its dual. Given a point $(p,H)$ of $A_{n,\{1,n\}}(\FF)$, let $\ell_p$ be the member of $\mathfrak{S}$ which contains $p$ and $L_H$ the member of $\mathfrak{S}^*$ contained in $H$. Then $H\supset \ell_p$ if and only if $p\in L_H$.
\end{lemma}

Under the hypotheses of Lemma \ref{main3 lem3}, we have 
\[\{(p,H)\in A_{n,\{1,n\}}(\FF) \mid H\supset \ell_p\} ~ = ~ \{(p,H)\in A_{n,\{1,n\}}(\FF) \mid p \in L_H\}\]
where $\ell_p$ and $L_H$ are defined as in the statement of Lemma \ref{main3 lem3} and we freely write $(p,H)\in A_{n,\{1,n\}}(\FF)$ to mean that $(p,H)$ is a point of $A_{n,\{1,n\}}(\FF)$. Put
\[\cH_{\mathfrak{S}} ~ := ~ \{(p,H)\in A_{n,\{1,n\}}(\FF) \mid H\supset \ell_p\} ~ (= \{(p,H)\in A_{n,\{1,n\}}(\FF) \mid p \in L_H\}).\]
The following will be proved in Section \ref{spread}.  

\begin{theo}\label{main3} 
Let $\mathfrak{S}$ be a line-spread of $\PG(n,\FF)$ and suppose that $\mathfrak{S}$ admits the dual. Then the set
$\cH_{\mathfrak{S}}$ defined as above is a hyperplane of $A_{n,\{1,n\}}(\FF)$.
\end{theo}  

We call the hyperplanes as $\cH_{\mathfrak{S}}$ {\em hyperlanes of spread type}. It is clear from the definition of $\cH_{\mathfrak{S}}$ that no hyperplane of spread type contains maximal singular subspaces of $A_{n,\{1,n\}}(\FF)$. So, in view of Corollary \ref{main2 bis} and Proposition \ref{hyp pure}, if $\cH_{\mathfrak{S}}$ arises from the natural embedding $\veS$, namely $\cH_{\mathfrak{S}} = \cH_M$ for a matrix $M\in M_{n+1}(\FF)\setminus[\lambda I]_{\lambda\in\FF}$, then $M$ admits no eigenvalue in $\FF$. More on $M$ when $\cH_M$ is of spread type will be stated below (Theorem \ref{main4}).

When $n = 3$ every line-spread which is also a dual line spread (as it is always the case when $\FF$ is finite) is its own dual. Let $n \geq  3$ be an odd integer and suppose that $\FF$ is not quadratically closed or it admits a non-trivial involutory automorphism. Then line-spreads of $\PG(n,\FF)$ exist which admit the dual (see Sections \ref{prel-spread-4} and \ref{prel-spread-3}; also below, Lemma \ref{main4 lem2}). In particular, every line-spread obtained by one of the constructions to be recalled in Sections \ref{prel-spread-1}, \ref{prel-spread-2} and \ref{prel-spread-3} admits the dual (see Section \ref{prel-spread-4}, Proposition \ref{dual cor} and Section \ref{prel-spread-3}, Proposition \ref{semilinear2}).     

We shall now turn to the last theorem of this section but, before to state it, we state the following two lemmas, to be proved in Section \ref{prel-spread-4}. 

\begin{lemma}\label{main4 lem1}
The following are equivalent for a matrix $M\in M_{n+1}(\FF)$:
\begin{itemize}
\item[$(\mathrm{S}_{\mathrm{mat}})$] $M^2\bx\in \langle \bx, M\bx\rangle$ for every non-zero vector $\bx\in V$;
\item[$(\mathrm{S}_{\mathrm{mat}}^*)$] $\xi M^2\in \langle \xi, \xi M\rangle$ for every non-zero linear functional $\xi\in V^*$.
\end{itemize}
\end{lemma}

\begin{lemma}\label{main4 lem2}
Let $M\in M_{n+1}(\FF)$ admit no eigenvalue in $\FF$ and suppose that $M$ satisfies the equivalent conditions $(\mathrm{S}_{\mathrm{mat}})$ and $(\mathrm{S}_{\mathrm{mat}}^*)$ of Lemma \ref{main4 lem1}. Put
\[\mathfrak{S}_M ~ := ~ \{\langle \bx, M\bx\rangle\}_{\bx\in V\setminus\{0\}}, \hspace{5 mm} \mathfrak{S}_M^* ~ := ~ \{\langle \xi, \xi M\rangle\}_{\xi\in V^*\setminus\{0\}}.\]
Then $\mathfrak{S}_M$ is line-spread of $\PG(n,\FF)$, it admits the dual and $\mathfrak{S}^*_M$ is its dual. 
\end{lemma} 

Lemmas \ref{main4 lem1} and \ref{main4 lem2} imply the following, to be proved in Section \ref{spread}.    

\begin{theo}\label{main4} 
A hyperplane $\cH_M$ of plain type is of spread type if and only if $M$ admits no eigenvalue in $\FF$ and satisfies conditions $(\mathrm{S}_{\mathrm{mat}})$ and $(\mathrm{S}_{\mathrm{mat}}^*)$ of Lemma \ref{main4 lem1}. If this is the case, then $\cH_M = \cH_{\mathfrak{S}_M}$, with $\mathfrak{S}_M$ as in Lemma \ref{main4 lem2}.   
\end{theo}

As it will be clear from Section \ref{prel-spread-2}, the line-spreads of $\PG(n,\FF)$ which can be obtained by the construction to be recalled in Section \ref{prel-spread-1}, called {\em spreads of standard type} in Section \ref{prel-spread}, are just the same as  the sets $\mathfrak{S}_M$ described in Lemma \ref{main4 lem2}, for a matrix $M$ admitting no eigenvalue in $\FF$ and satisfying condition $(\mathrm{S}_{\mathrm{mat}})$ of Lemma \ref{main4 lem1}.  
Accordingly, the hyperplanes of spread type which arise from the natural embedding $\veS$ are precisely those that can be obtained from spreads of standard type. However, not every line-spread which amits the dual is of standard type (see Section \ref{prel-spread-3}, Proposition \ref{semilinear2}). So, hyperplanes of spread type also exist which do not arise from $\veS$. 

\paragraph{Organization of the paper.} In Section \ref{prel}, after having fixed some notation,  we discuss a few properties of the saturation form $f$, some of which have been already mentioned in Section \ref{tensor}. Next we turn to line-spreads of projective spaces, recalling some information on this topic. Section \ref{pure} contains the proofs of Propositions \ref{hyp pure} and \ref{main1}, Section \ref{max} the proof of Theorem \ref{main2} and Section \ref{spread} the proofs of Theorems \ref{main3} and \ref{main4}. 

\section{Preliminaries}\label{prel}

\subsection{Notation} 

\subsubsection{Notation for vectors, tensors and matrices}\label{notation1}

So far we have have implicitly adopted the convention of denoting the vectors of $V = V(n+1,\FF)$ by low case boldface roman letters, the vectors of its dual $V^*$ by low case greek letters and the matrices of $M_{n+1}(\FF)$, namely the vectors of $V\otimes V^*$, by capital roman letters, but we used $0$ to denoted both the null-vector of $V$ and the null-vector of $V^*$. We shall keep these conventions througout this paper. Scalars will be denoted by roman or greek low case letters. As in the previous section, $O$ and $I$ are the null matrix and the identity matrix of $M_{n+1}(\FF)$. 

Given a basis $E = (\be_i)_{i=0}^n$ of $V$ let $E^* = (\eta_i)_{i=0}^n$ be the basis of $V^*$ dual to it. So, $\eta_i(\be_j) = \delta_{i,j}$ (Kronecker symbol) for any choice of $i, j = 0, 1,..., n$. The pure tensors $E_{i,j} = \be_i\otimes \eta_j$ form a basis of $V\otimes V^*$, henceforth denoted by $E\otimes E^*$. When $E$ is the natural basis of $V = V(n+1,\FF)$ then $E\otimes E^*$ is the natural basis of $V\otimes V^*$ and $E_{i,j}$, regarded as a matrix, is the matrix where all entries are null but the $(i,j)$-entry, which is $1$. So, $E\otimes E^*$ is just the usual natural basis of $M_{n+1}(\FF)$.    

We denote by $\Pu(V\otimes V^*)$ ($= \Pu(M_{n+1}(\FF))$) the set of pure tensors of $V\otimes V^* = M_{n+1}(\FF)$, with the convention that the null vector $O$ of $V\otimes V^*$ is not included in that set. So, $\Pu(M_{n+1}(\FF))$ is the set of matrices of $M_{n+1}(\FF)$ of rank 1. The hyperplane $M^0_{n+1}(\FF)$ of $M_{n+1}(\FF)$ will also be denoted by $(V\otimes V^*)_0$ and we set $\Pu((V\otimes V^*)_0) := \Pu(V\otimes V^*)\cap(V\otimes V^*)_0$. So, $\Pu(M_{n+1}^0(\FF))$ ($= \Pu((V\otimes V^*)_0)$) is the set of null-traced matrices of $M_{n+1}(\FF)$ of rank $1$, namely the set of pure tensors $\bx\otimes \xi$ with $\xi(\bx) = 0$.

Given two non-zero vectors $\bx, \by \in V$, if $\bx$ and $\by$ are proportional we write $\bx \equiv \by$. We use the same notation for vectors of $V^*$ and matrices of $M_{n+1}(\FF)$. For instance, for a matrix $M\in M_{n+1}(\FF)\setminus\{O\}$, when writing $M \equiv I$ we mean that $M$ is a scalar matrix. 

When we need to distinguish betwen a non-zero vector $\bx$ of $V$ and the point of $\PG(V) = \PG(n,\FF)$ represented by it, we denote the latter by $[\bx]$. Similarly, for a set $X$ of vectors of $V$, we put 
\begin{equation}\label{[X]}
[X] ~ := ~ \{[\bx] \mid \bx\in X\setminus\{0\}\}.
\end{equation} 
In particular, if $X$ is a subspace of $V$ then $[X]$ is the subspace of $\PG(V)$ corresponding to it. The same conventions will be adopted for $V^*$ and $V\otimes V^*$. 

Finally, we state some terminology and notation for collineations of $\PG(V)$. Given a semi-linear invertible transformation $g:V\rightarrow V$, we denote by $[g]$ the collineation of $\PG(V)$ corresponding to $g$. We say that a collineation $\gamma$ of $\PG(V)$ is {\em linear} if $\gamma = [g]$ for an invertible linear trasformation $g:V\rightarrow V$.  

Given a linear transformation $g:V\rightarrow V$ and a basis $E = (\be_i)_{i=0}^n$ of $V$, let $M = (m_{i,j})_{i,j=0}^n$ be the representative matrix of $g$ with respect to $E$. So, for a vector $\bx = \sum_{i=0}^n\be_ix_i$ of $V$, we have $f(\bx) = \sum_{i,j =0}^n\be_im_{i,j}x_j$. Regarded the $(n+1)$-tuples $(\sum_{j=0}^nm_{i,j}x_j)_{i=0}^n$ and $(x_i)_{i=0}^n$ as $(n+1)\times 1$ matrices, we have $(\sum_{j=0}^nm_{i,j}x_j)_{i=0}^n = M\cdot(x_i)_{i=0}^n$. We keep a record of this fact by writing $g(\bx) = M\bx$, as if $E$ was the natural basis of $V$, also replacing the symbol $g(\bx)$ with $M\bx$ whenver this matrix-style notation is convenient. We adopt the same convention for linear mappings $g^*:V^*\rightarrow V^*$ with the dual $E^*$ of $E$ as a basis of $V^*$, but for regarding the vectors of $V^*$ as $(1\times(n+1))$-matrices, thus writing $\xi M$ for $g^*(\xi)$. So, the same matrix $M\in M_{n+1}(\FF)$ represents both an endomorphism $g$ of $V$ (with respect to $E$) and an endomorphism $g^*$ of $V^*$ (with respectto $E^*$). 

\subsubsection{Notation for the geometry $A_{n,\{1,n\}}(\FF)$}\label{notation2} 

Henceforth the symbol $\perp$ (not to be confused with the symbol $\perp_f$ introduced in Section \ref{tensor}) stands for collinearity in $A_{n,\{1,n\}}(\FF)$ and we adopt for $\perp$ all usual conventions from the theory of point-line geometries. Thus, given a point $(a,A)$ of $A_{n,\{1,n\}}(\FF)$, we denote by $(a, A)^\perp$ the set of points of $A_{n,\{1,n\}}(\FF)$ collinear with $(a,A)$, the point $(a,A)$ being included among them by convention. Given a set $\cal X$ of points of $A_{n,\{1,n\}}(\FF)$ we put ${\cal X}^\perp := \cap_{(p,H)\in{\cal X}}(p,H)^\perp$. 

With $\cal X$ as above, we denote by $\langle {\cal X}\rangle$ the subspace of $A_{n,\{1,n\}}(\FF)$ generated by $\cal X$, namely the smallest subspace of $A_{n,\{1,n\}}(\FF)$ which contains $\cal X$. We warn that the symbol $\langle . \rangle$ will also be used to denote spans in $V$, $V^*$, $V\otimes V^*$ and in their projective spaces $\PG(V)$, $\PG(V^*)$ and $\PG(V\otimes V^*)$. However, in any case the context will make it clear which kind of span we are referring to. 

\subsection{Properties of the saturation form}\label{prel-form}

Let $f$ be the saturation form of $M_{n+1}(\FF)$, as defined in (\ref{f1}) and (\ref{f2}) of Section \ref{tensor}. Formula (\ref{f2}) shows that, with respect to a suitable ordering of the natural basis of $M_{n+1}(\FF)$, the form $f$ is represeted by a block-diagonal matrix, with $(n+1)n/2$ blocks of order $2$ equal to
\[\left(\begin{array}{cc}
0 & 1 \\
1 & 0 
\end{array}\right)\]
and one more block equal to the identity matrix $I$ of order $n+1$. This matrix is non-singular. Hence $f$ is non-degenerate, as claimed in Section \ref{tensor}. We denote by $\cS_f$ the set of $f$-isotropic vectors of $M_{n+1}(\FF)$, namely the set of matrices $X\in M_{n+1}(\FF)$ such that $f(X,X) = 0$.  

If $\mathrm{rank}(X) = 1$, namely $X = \bx\otimes \xi \in \Pu(V\otimes V^*)$, then $f(X,X) = (\xi(\bx))^2$ by formula (\ref{f3}) of Section \ref{tensor}.  Hence $f(X,X)= 0$ if and only if $\xi(\bx) = 0$, namely $X\in M_{n+1}^0(\FF)$. Accordingly,  
\[\cS_f\cap \Pu(V\otimes V^*) ~ = ~ \Pu((V\otimes V^*)_0).\] 
Note that $[\Pu((V\otimes V^*)_0)]$ is just the $\veS$-image of the point-set of $A_{n,\{1,n\}}(\FF)$. 

\begin{note}\label{quadric}
\em 
When $\mathrm{char}(\FF) \neq 2$ the set $\cS_f$ spans $M_{n+1}(\FF)$. In this case the points of $\PG(M_{n+1}(\FF))$ represented by the non-zero vectors of $\cS_f$ form a non-singular quadric. On the other hand, let $\mathrm{char}(\FF) = 2$. Then $f(X,X) = 0$ if and only if $\mathrm{trace}(X) = 0$.  So, if $\mathrm{char}(\FF) = 2$ then $\cS_f = M_{n+1}^0(\FF)$. With the terminology of Tits \cite[chp. 8]{T}, the form $f$ is trace-valued or non-trace-valued according to whether $\mathrm{char}(\FF) \neq 2$ or $\mathrm{char}(\FF) = 2$, but this discrepancy has no relevance for the matter of this paper. 
\end{note}

As in Section \ref{tensor}, we denote by $\perp_f$ the orthogonality relation associated to $f$. As $f$ is non-degenerate, the hyperplanes of $M_{n+1}(\FF)$ are the perps $M^{\perp_f}$ for $M\in M_{n+1}(\FF)\setminus\{O\}$ and, for two matrices $M, N\in M_{n+1}(\FF)\setminus\{O\}$, we have $M^{\perp_f} = N^{\perp_f}$ if and only if $M \equiv N$. It is clear from formula (\ref{f1}) of Section \ref{tensor} that   $I^{\perp_f} = M_{n+1}^0(\FF)$. Therefore, for $M\in M_{n+1}(\FF)\setminus\{O\}$, we have $M^{\perp_f} = M_{n+1}^0(\FF)$ if and only if $M\equiv I$.

Every hyperplane of $M_{n+1}^0(\FF)$ is the intersection of $M^0_{n+1}(\FF)$ with a hyperplane $M^{\perp_f}$ of $M_{n+1}(\FF)$ for a suitable matrix $M \not\in \langle I\rangle$. Of course the equality $M^{\perp_f}\cap M_{n+1}^0(\FF) = N^{\perp_f}\cap M_{n+1}^0(\FF)$ does not imply that $M\equiv N$. The next statement, already mentioned in Remark \ref{rem notation}, follows from well known properties of the polarities associated to non-degenerate reflexive sesqulinear forms.  

\begin{prop}\label{prop0}
For $M, N\in M_{n+1}(\FF)\setminus\{O\}$, we have $M^{\perp_f}\cap M_{n+1}^0(\FF) = N^{\perp_f}\cap M_{n+1}^0(\FF)$ if and only if 
$\langle M, I\rangle = \langle N, I\rangle$.  
\end{prop}

\subsection{Line-spreads}\label{prel-spread}

In the first part of this subsection (Sections \ref{prel-spread-1} and \ref{prel-spread-2}) we survey a few basics on line-spreads, focusing on a particular but important class of them. All we are going to say in this part is well known (see e.g. \cite{J}). Nevertheless, for most of the facts we are going to mention in our survey we will give at least a sketch of the proof, in order to make our exposition as self contained as possible. In Section \ref{prel-spread-4} we turn to duals of line-spreads, proving a number of claims made in Section \ref{intro spread}. Finally, in Section \ref{prel-spread-3} we turn back to a construction of spreads from suitable linear collineations, discussed in Section \ref{prel-spread-2}, generalzing it for non-linear collineations.          

Throughout this subsection $n > 2$ is odd. In Sections \ref{prel-spread-1}, \ref{prel-spread-2} and \ref{prel-spread-4} we also assume that $\FF$ is not quadratically closed. This hypothesis is dropped in Section \ref{prel-spread-3}.    

\subsubsection{A standard way to construct line-spreads}\label{prel-spread-1}

Let $\ovF$ be an extension of $\FF$ of degree $2$ and  $V = V(n+1, \FF)$, as before. As $n$ is odd, the ratio $m := (n+1)/2$ is an integer.  Put $\ovV = V(m, \ovF)$. For every non-zero vector $\bar{\bv}\in\overline{V}$, the $1$-dimensional subspace $\langle \bar{\bv}\rangle$ of $\ovV$ can be regarded as a $2$-dimensional subspace over $\FF$. Given any two scalars $a, b \in \overline{\FF}$, non-proportional as vectors of the $\FF$-space $\ovF$, the vectors $\bar{\bv}a$ and $\bar{\bv}b$ form a basis of the $\FF$-space $\langle\bar{\bv}\rangle$. Accordingly, chosen an ordered basis $(\bar{\be}_j)_{j=1}^m$ of $\ovV$ (its natural basis, for instance), we can regard $\overline{V}$ as an ($n+1)$-dimensional $\FF$-space with $(\bar{\be}_1a_1, \bar{\be}_1b_1,..., \bar{\be}_ma_m, \bar{\be}_mb_m)$ as a basis, where for every $j = 1,..., m$ the scalars $a_j$ e $b_j$ are non-proportional as vectors of the $\FF$-space $\overline{\FF}$. With $a_1, ..., a_m, b_1,..., b_m$ as above, we say that the $m$-tuples $\mathfrak{a} = (a_j)_{j=1}^m$ and $\mathfrak{b} = (b_j)_{j=1}^m$ are {\em elementwise non-proportional} ({\em over $\FF$}).

Given a basis $E = (\be_i)_{i=0}^n$ of $V$ and two elementwise non-proportional $m$-tuples $\mathfrak{a} = (a_j)_{j=1}^m$ and $\mathfrak{b} = (b_j)_{j=1}^m$, a unique isomorphism from $V$ and the $\FF$-space $\overline{V}$ exists which for $j = 1,..., m$ maps $\be_{2j-2}$ onto $\bar{\be}_ja_j$ and $\be_{2j-1}$ onto $\bar{\be}_{j}b_j$. We denote this isomorphism by $\phi_{E, \mathfrak{a},\mathfrak{b}}$.  

With $\phi_{E,\mathfrak{a},\mathfrak{b}}$ as above, for every non-zero vector $\bar{\bv}\in \overline{V}$, the $\phi_{E,\mathfrak{a},\mathfrak{b}}$-preimage $\phi_{E,\mathfrak{a},\mathfrak{b}}^{-1}(\langle \bar{\bv}\rangle)$ of the $1$-dimensional subspace of the $\overline{\FF}$-space $\overline{V}$ spanned by $\bar{\bv}$ is a $2$-dimensional subspace of $V$, hence a line of $\PG(V)$. The set of lines of $\PG(V)$ obtained in this way is a spread of $\PG(V)$, which we shall denote by $\mathfrak{S}_{E,\mathfrak{a},\mathfrak{b}}$. We call the line-spreads constructed in this way {\em spreads of standard type defined by means of} $\ovF$, also {\em standard spreads} for short, when the particular choice of the extension $\ovF$ is clear from the context or irrelevant.    

In particular, let $\mathfrak{u} = (1,1...,1)$ be the $m$-tuple with all terms equal to $1$ and, given $\omega\in \ovF\setminus\FF$, let $\mathfrak{u}\omega := (\omega, \omega,..., \omega)$. Clearly, $\mathfrak{u}$ and $\mathfrak{u}\omega$ are elementwise non-proportional. So,  $\mathfrak{S}_{E, \mathfrak{u}, \mathfrak{u}\omega}$ is a spread of standard type, henceforth denoted by $\mathfrak{S}_{E,\omega}$ (also $\mathfrak{S}_\omega$ for short). 

We say that two line-spreads $\mathfrak{S}$ and $\mathfrak{S}'$ of $\PG(V)$ are {\em equivalent} if $\PG(V)$ admits a linear collineation $\gamma$ such that $\gamma(\mathfrak{S}) = \mathfrak{S}'$. If $\mathfrak{S}$ and $\mathfrak{S}'$ are equivalent we write $\mathfrak{S}\sim \mathfrak{S}'$.  

\begin{prop}\label{smart 00}
All line-spreads of $\PG(V)$ of standard type defined by means of the same extension $\ovF$ of $\FF$ are mutually equivalent. In particular, chosen a scalar $\omega\in \ovF\setminus\FF$ and a basis $E$ of $V$, all of them are equivalent to $\mathfrak{S}_{E,\omega}$.   
\end{prop}
{\bf Proof.} Given two bases $E = (\be_i)_{i=0}^n$ and $E'$ of $V$, let $\mathfrak{S} = \mathfrak{S}_{E',\mathfrak{a},\mathfrak{b}}$ for a pair of elementwise non-proportional $m$-tuples $\mathfrak{a} = (a_j)_{j=1}^m$ and $\mathfrak{b} = (b_j)_{j=1}^m$ and let $\omega\in\ovF\setminus\FF$. As $V$ admits an automorphism which maps $E'$ onto $E$, we can assume with no loss that $E' = E$. The $\ovF$-space $\ovV$ admits and $\FF$-linear mapping $\bar{g}$ which maps $(\bar{e}_ja_j, \bar{e}_jb_j)$ onto $(\bar{e}_j, \bar{e}_j\omega)$, for every $j = 1,..., m$. The composite mapping $g := \phi_{E,\mathfrak{u}, \mathfrak{u}\omega}^{-1}\cdot \bar{g}\cdot\phi_{E,\mathfrak{a},\mathfrak{b}}$ maps $\mathfrak{S}$ onto $\mathfrak{S}_{E,\omega}$.  \hfill $\Box$ 

\subsubsection{Constructing line-spreads by means of collineations}\label{prel-spread-2}

As in the previous subsection, $\ovF$ is an extension of $\FF$ of degree $2$. We recall that the {\em characteristic polynomial} of a scalar $\omega\in\ovF\setminus\FF$ is the unique monic polynomial $P_\omega(t) = t^2+at+b$ such that $P_\omega(\omega) = 0$ but the equation $t^2+at+b = 0$ has no solution in $\FF$.   

\begin{lemma}\label{smart 010}
Given a scalar $\omega\in \ovF\setminus\FF$ and a basis $E = (\be_i)_{i=1}^n$ of $V$, there exists an invertible linear mapping $f_\omega:V\rightarrow V$  such that all the following hold. 
\begin{itemize}
\item[$(1)$] $\mathfrak{S}_{E,\omega} = \{\langle [\bv], [f_\omega(\bv)]\rangle\}_{\bv\in V\setminus\{0\}}$.
\item[$(2)$] Let $M_\omega\in M_{n+1}(\FF)$ be the matrix which represents $f_\omega$ with respect to $E$. Then $M_\omega$ is block-diagonal with $m$ blocks $A_1,..., A_m$ of order $2$, all of wich are equal to 
\[\left(\begin{array}{cr} 
0 & -b\\
1 & -a
\end{array}\right)\]
where $t^2+at+b = P_\omega(t)$. 
\item[$(3)$] The minimal and characteristic polynomials of $f_\omega$ are equal to $P_\omega(t)$ and $(P_\omega(t))^m$ respectively.  Accordingly, $f_\omega$ admits no eigenvalue in $\FF$.
\item[$(4)$] $f^2_\omega(\bv) \in \langle \bv, f_\omega(\bv)\rangle$ for every $\bv\in V$. Equivalently, the collineation $[f_\omega]$ of $\PG(V)$ represented by $f_\omega$ stabilizes all lines of $\mathfrak{S}_{E,\omega}$.  
\end{itemize}
\end{lemma}
{\bf Proof.} Let $\bar{f}_\omega:\ovV\rightarrow\ovV$ be the $\ovF$-linear mapping which maps every vector $\bar{v}\in \ovV$ onto $\bar{v}\omega$. The $\FF$-linear mapping $f_\omega := \phi^{-1}_{E, \mathfrak{u}, \mathfrak{u}\omega}\cdot\bar{f}_\omega\cdot \phi_{E, \mathfrak{u}, \mathfrak{u}\omega}$ has the required properties. \hfill $\Box$ 

\begin{prop}\label{smart 01}
For every line-spread $\mathfrak{S}$ of $\PG(V)$ of standard type there exists an invertible linear mapping $f:V\rightarrow V$ such that:
\begin{itemize}
\item[$(1)$] $\mathfrak{S} = \{\langle [\bv], [f(\bv)]\rangle\}_{\bv\in V\setminus\{0\}}$.
\item[$(2)$]  $f$ admits no eigenvalue in $\FF$.
\item[$(3)$] $f^2(\bv) \in \langle \bv, f(\bv)\rangle$ for every $\bv\in V$. 
\end{itemize} 
\end{prop}
{\bf Proof.} In view of Proposition \ref{smart 00}, given $\omega\in\ovF\setminus\FF$, an invertible linear mapping $g:V\rightarrow V$ exists such that $[g]$ maps $\mathfrak{S}$ onto $\mathfrak{S}_\omega$. With $f_\omega$ as in Lemma \ref{smart 010}, the linear mapping $f := g^{-1}f_\omega g$ has the required properties. \hfill $\Box$   

\begin{note}\label{non-uniqueness}
\em
The mappings satisfying properties (1), (2) and (3) of Proposition \ref{smart 01} are those which represent collineations of $\PG(V)$ which stabilize $\mathfrak{S}$ line-wise and act fixed-point-freely on $\PG(V)$. For a given spread $\mathfrak{S}$ different collineations exist with these properties. Explicitly, if two linear mappngs $f$ and $f'$ of $V$ define line-spreads of $\PG(V)$, they define the same spread if and only if $f-f' \in \langle \mathrm{id}_V\rangle$.        
\end{note}

Conversely, let $\varphi$ be a linear collineation of $\PG(V)$ and suppose that $\varphi$ acts fixed-point freely on $\PG(V)$. Equivalently, the linear mappings which represent $\varphi$ admit no eigenvalue in $\FF$. Therefore $\ell_{p,\varphi} := \langle p, \varphi(p)\rangle$ is a line of $\PG(V)$, for every point $p\in \PG(V)$. Put
\begin{equation}\label{S-phi}
\mathfrak{S}_\varphi ~ := ~ \{\ell_{p,\varphi}\}_{p\in \PG(V)}.
\end{equation} 
Clearly, $\mathfrak{S}_\varphi$ covers the set of points of $\PG(V)$ but in general it is not a line-spread of $\PG(V)$. The followig is obvious.  

\begin{prop}\label{smart 1}
The set $\mathfrak{S}_\varphi$ is a line-spread if and only if 
\begin{equation}\label{eq-phi}
\varphi^2(p) ~\in ~\langle p, \varphi(p)\rangle ~  \mbox{ for every point } p\in \PG(V).
\end{equation}  
\end{prop}

So, given a linear collineation $\varphi$ of $\PG(V)$ acting fixed-point freely on $\PG(V)$ and satisfying (\ref{eq-phi}), the set of lines $\mathfrak{S}_\varphi$ defined as in (\ref{S-phi}) is a line-spread of $\PG(V)$. However, this construction yields no new spreads. Indeed let $f_\varphi$ be a representative of $\varphi$ in the group of linear automorphisms of $V$. Then we can always choose a basis $E$ of $V$ such that the matrix which represents $f_\varphi$ with respect to $E$ is as described in $(2)$ of Lemma \ref{smart 010}, for a suitable equation $t^2+at+b =0$ with no solutions in $\FF$. Of course, $\ovF := \FF(\omega)$ with $\omega^2 = -a\omega -b$ is the extension of $\FF$ which suits $\mathfrak{S}_\varphi$ and $\mathfrak{S}_\varphi =  \mathfrak{S}_{E,\omega}$. We have proved the following.  

\begin{prop}\label{smart 2}
With $\varphi$ as in Proposition \ref{smart 1}, the line-spread $\mathfrak{S}_\varphi$ is of standard type. 
\end{prop}

\begin{cor}
Let $\varphi$ and $\psi$ be two linear collineations of $\PG(V)$ acting fixed-point freely on $\PG(V)$ and satisfying {\rm (\ref{eq-phi})} of Proposition \ref{smart 1}. Let  $f_\varphi$ and $f_\psi$ be linear mappings of $V$ representing $\varphi$ and $\psi$ respectively and let $P_\varphi$ and $P_\psi$ be their minimal  polynomials (which have degree $2$ and admit no zeros in $\FF$). Then $\mathfrak{S}_\varphi\sim \mathfrak{S}_\psi$ if and only if $P_\varphi$ and $P_\psi$ define the same extension of $\FF$.
\end{cor}

\subsubsection{The dual of a line-spread}\label{prel-spread-4}  

As recalled in Remark \ref{composition}, a line-spread admits the dual only if it is a composition line-spread. \textcolor{red}{Composition line-spreads are well known objects, thoroughly studied by a number of authors.} Some of the proofs we are going to give in this subsection can be shortened or even dropped by referring to properties of composition spreads. However, since all of these proofs are easy, we prefer to expose them in full length, without asking the reader to look for appropriate references in the literature.

Our definition ofthe dual of a line spread is based on two properties, called $(\mathrm{S})$ and $(\mathrm{S}^*)$ in Section \ref{intro spread}. In Lemma \ref{main3 lem1} we claimed that these two properties are equivalent but we didn't prove this claim. We prove it here.\\

\noindent
{\bf Proof of Lemma \ref{main3 lem1}.} Suppose that $\mathfrak{S}$ and $\mathfrak{S}^*$ satisfy $(\mathrm{S})$. We shall prove that they also satisfy $(\mathrm{S}^*)$. Let $\ell\in\mathfrak{S}$, let $H$ be a hyperplane of $\PG(V)$ and $L$ the unique sub-hyperplane of $\mathfrak{S}^*$ contained in $H$. According to $(\mathrm{S}^*)$, we must have $\ell \subseteq L$. This is indeed the case. Indeed $H$ contains both $\ell$ and $L$. Hence $\ell\cap L\neq \emptyset$. Let $p\in \ell\cap L$. According to $(\mathrm{S})$, the sub-hyperplane $L$ contains the unique line of $\mathfrak{S}$ through $p$. However that line is just $\ell$. Hence $L\supseteq \ell$. So, $(\mathrm{S})$ does imply $(\mathrm{S}^*)$. The converse implication, from $(\mathrm{S}^*)$ to $(\mathrm{S})$, can be proved by dualizing the above argument.  \hfill $\Box$ \\

In Lemma \ref{main3 lem2} we claim that every line-spread admits at most one dual. Here is a proof of that claim.\\

\noindent
{\bf Proof of Lemma \ref{main3 lem2}.} For a contradiction, suppose that two disctinct dual line-spreads $\mathfrak{S}^*_1$ and $\mathfrak{S}^*_2$ exist such that property $(\mathrm{S})$ holds for both pairs $\{\mathfrak{S}, \mathfrak{S}^*_1\}$ and $\{\mathfrak{S}, \mathfrak{S}^*_2\}$. As $\mathfrak{S}^*_1\neq \mathfrak{S}^*_2$, there exists a hyperplane $H$ of $\PG(n,\FF)$ such, if $L_i$ is the member of $\mathfrak{S}^*_i$ contained in $H$ for $i = 1, 2$, then $L_1\neq L_2$. Choose $p \in H\setminus(L_1\cap L_2)$ and let $\ell_p$ be the line of $\mathfrak{S}$ through $p$. Then $\ell_p\cap L_i\neq \emptyset$ for $i =1, 2$, as $L_i$ is a hyperlane of $H$. Property $(\mathrm{S})$ then forces $\ell_p  \subseteq L_i$. Hence $\ell_p\subseteq L_1\cap L_2$ and therefore $p\in L_1\cap L_2$. This contradicts the choice of $p$. \hfill $\Box$\\

\noindent
{\bf Proof of Lemma \ref{main3 lem3}.} Let $\mathfrak{S}$ admit the dual and let $\mathfrak{S}^*$ be its dual. For a point-hyperplane flag $(p,H)$ of $\PG(\FF)$, suppose that $H$ contains the line $\ell_p\in\mathfrak{S}$ through $p$. By $(\mathrm{S}^*)$, the sub-hyperplane $L_H\in\mathfrak{S}^*$ contained in $H$ contains $\ell_p$. Hence $p\in L_H$. Conversely, if $p\in L_H$ then $\ell_p\subseteq L_H$ by $(\mathrm{S})$ and therefore $H\supset \ell_p$. \hfill $\Box$ \\

Two more lemmas stated in Section \ref{intro spread} remain to be proved, namely Lemmas \ref{main4 lem1} and \ref{main4 lem2}.\\

\noindent
{\bf Proof of Lemma \ref{main4 lem1}.} Suppose that $M^2\bx \in \langle \bx, M\bx\rangle$ for every $\bx\in V$. Then, given $\xi\in V^*$, for every $\bx\in V$, if $\xi(\bx) = \xi(M\bx) = 0$ then $\xi(M^2\bx) = 0$. However $\xi(M\bx) = (\xi M)(\bx)$ and $\xi(M^2\bx) = (\xi M^2)(\bx)$. So, the previous implication can also be read as follows: $(\xi M^2)(\bx) = 0$ whenever $\xi(\bx) = (\xi M)(\bx) = 0$, which holds true if and only if $\xi M^2\in \langle \xi, \xi M\rangle$. This proves that claim $(1)$ of Lemma \ref{main4 lem1} implies $(2)$ of that lemma. The converse implication from $(2)$ to $(1)$ is proved by reversing the previous argument.   \hfill $\Box$\\

\noindent
{\bf Proof of Lemma \ref{main4 lem2}.} Let $f_M:V\rightarrow V$ be the linear mapping represented by $M$ with respect to the natural basis $E$ of $V$ and $\varphi$ the corresponding collineation of $\PG(V)$. The hypotheses assumed on $M$ in Lemma \ref{main4 lem2} amount to say that $f_M$ admits no eigenvalues in $\FF$ and $f_M^2(\bx) \in\langle \bx, f_M(\bx)\rangle$ for every $\bx\in V$. So, the set of lines $\mathfrak{S}_M$ is the same as the line-spread of standard type $\mathfrak{S}_\varphi$ associated to the collineation $\varphi = [f_M]$ (compare Propositions \ref{smart 01} and \ref{smart 1}). Similarly, $\mathfrak{S}^*_M$ is the standard spread $\mathfrak{S}_{\varphi^*}$ of $\PG(V^*)$, where $\varphi^*$ is the collineation of $\PG(V^*)$ corresponding to the linear mapping $f_M^*:V^*\rightarrow V^*$ represented by $M$ with respect to the basis $E^*$ of $V^*$ dual of the given basis $E$ of $V$. So, $\mathfrak{S}_M$ and $\mathfrak{S}^*_M$ are indeed line-spreads of $\PG(V)$ and $\PG(V^*)$ respectively, as claimed in Lemma \ref{main4 lem2}. In fact, they are spreads of standard type. 

It remains to show that the pair $\{\mathfrak{S}_M, \mathfrak{S}^*_M\}$ satisfies property $(\mathfrak{S})$. In view of the way the lines of $\mathfrak{S}_M$ and the dual lines of $\mathfrak{S}^*_M$ are defined, this amounts to the following: for every $\bx\in V\setminus\{0\}$ and every $\xi\in V\setminus\{0\}$, if $\xi(\bx) = (\xi M)\bx = 0$ then $\xi(M\bx) = (\xi M)(M\bx)$. In other words, since $(\xi M)\bx = \xi(M\bx)$ and $(\xi M)(M\bx) = \xi(M^2\bx)$, if $\xi(\bx) = \xi(M\bx) = 0$ then $\xi(M^2\bx) = 0$. This implication follows from the hypothesis that $M^2\bx\in \langle \bx, M\bx\rangle$. \hfill $\Box$ \\   

In view of Proposition \ref{smart 01}, Lemma \ref{main4 lem2} implies the following.

\begin{prop}\label{dual cor}
Every line-spread of standard type admits the dual.
\end{prop}  

However, spreads of non-standard type also exist which admit the dual. We shall describe a class of them in the next subsection. 

\subsubsection{Line-spreads of semi-standard type and their duals}\label{prel-spread-3}

In Proposition \ref{smart 1} we claim that, for a linear collineation $\varphi$ of $\PG(V)$ acting fixed-point freely on $\PG(V)$, the set of lines $\mathfrak{S}_\varphi$ defined as in (\ref{S-phi}) is a line-spread if and only if $\varphi$ satisfies (\ref{eq-phi}). 
However, keeping the hypothesis that $\varphi$ acts fixed-point freely on $\PG(V)$, the statement of Proposition \ref{smart 1} holds true even if $\varphi$ is non-linear (and even if the field $\FF$ is quadratically closed). We say that a line-spread $\mathfrak{S}_\varphi$ with $\varphi$ a non-linear collineation of $\PG(V)$ is a spread of {\em semi-standard type}. 

\begin{lemma}\label{semilinear1}
A non-linear collineation of $\PG(V)$ acting fixed-point freely on $\PG(V)$ defines a line-spread of $\PG(V)$ if and only if it is an involution.
\end{lemma} 
{\bf Proof.} The `if' part of the lemma is obvious. Turning to the `only if' part, let $\varphi$ be a non-linear collineation of $\PG(V)$ which fixes no points of 
$\PG(V)$ and suppose that $\mathfrak{S}_\varphi$ is a spread. Let $\sigma$ be the (non-trivial) automorphism of $\FF$ assocated to $\varphi$ and $f:V\rightarrow V$ a semi-linear transformation of $V$ corresponding to $\varphi$. As $\mathfrak{S}_\varphi$ is a spread, for every vector $\bx\in V\setminus\{0\}$ there exist two scalars $\lambda_\bx, \mu_\bx\in \FF$ such that 
\begin{equation}\label{semi 1}
f^2(\bx) ~ = ~ \bx\lambda_\bx + f(\bx)\mu_\bx
\end{equation}
and at least one of them is different from $0$. Equation (\ref{semi 1}) yields 
\begin{equation}\label{semi 2}
\begin{array}{rcl}
f^2(\bx+\by) & = & f^2(\bx) + f^2(\by) =  \bx\lambda_\bx + \by\lambda_\by +  f(\bx)\mu_\bx + f(\by)\mu_\by,\\
f^2(\bx+\by) & = & (\bx+\by)\lambda_{\bx+\by} + (f(\bx)+f(\by))\mu_{\bx+\by}.
\end{array} 
\end{equation}
Suppose that $\by\not\in\langle \bx, f(\bx)\rangle$. Since $\mathfrak{S}_\varphi$ is a spread, the set $\{\bx, f(\bx), \by, f(\by)\}$ is linearly independent. Hence the two equations of (\ref{semi 2}) imply that $\lambda_\bx = \lambda_{\bx+\by} =\lambda_\by$ and $\mu_\bx = \mu_{\bx+by} = \mu_\by$. Consequently, the scalars $\lambda_\bx$ and $\mu_\bx$ in (\ref{semi 1}) do not depend on the particular choice of the vector $\bx$ and (\ref{semi 1}) can be written as follows:
\begin{equation}\label{semi 3}
f^2(\bx) ~ = ~ \bx\lambda + f(\bx)\mu ~ \mbox{ for every } \bx\in V. 
\end{equation}
By (\ref{semi 3}), the following holds for every $t\in \FF$:
\begin{equation}\label{semi 4}
\bx\lambda t^{\sigma^2} + f(\bx)\mu t^{\sigma^2} ~ =~ f^2(\bx)t^{\sigma^2} ~= ~f^2(\bx t) ~ = ~ \bx t\lambda  + f(\bx)t^\sigma\mu. 
\end{equation} 
If $\mu \neq 0$ equation (\ref{semi 4}) implies that $t^{\sigma^2} = t^\sigma$ for every $t\in \FF$. This cannot be, since $\sigma\neq \mathrm{id}_\FF$.
Therefore $\mu = 0$ and (\ref{semi 3}) just says that $\phi^2$ is the identity.  \hfill $\Box$  

\begin{prop}\label{semilinear2}
Every semi-standard line-spread admits the dual.
\end{prop}
{\bf Proof.} Let $\varphi$ be a non-linear collineation of $\PG(V)$. So, $\varphi = [f]$ for a semi-linear transformation $f:\bx\rightarrow M\bx^\sigma$ of $V$, where $M\in M_{n+1}(\FF)$ is non-singular, $\sigma$ is the non-trivial automorphism of $\FF$ assciated to $\varphi$ and, for $\bx = (x_i)_{i=0}^n \in V$ we put $\bx^\sigma = (x_i^\sigma)_{i=0}^n$. Suppose that $\varphi$ defines a line-spread $\mathfrak{S}_\varphi$. Then $f^2(\bx) \equiv \bx$ for every $\bx\in V$ by Lemma \ref{semilinear1}. Therefore 
\begin{equation}\label{inv2}
\sigma^2 = \mathrm{id}_\FF ~ \mbox{ and } ~ M^\sigma \equiv M^{-1}.
\end{equation} 
We shall prove that the following also holds true. 	
\begin{equation}\label{inv3}
\xi M \not\equiv \xi^\sigma  ~ \mbox{ for every } \xi\in V^*\setminus\{0\}, 
\end{equation}
By way of contradiction, suppose that a linear functional $\alpha\in V^*\setminus\{0\}$ exists such that $\alpha M \equiv\alpha^\sigma$. Then the conditions $\alpha M\bx^\sigma = 0$ and $\alpha\bx = 0$ are equivalent for every $\bx\in V$. However, if $A$ is the hyperplane of $\PG(V)$ represented by $\alpha$ and $\bx\neq 0$, the conditions $\alpha M\bx^\sigma = 0$ and $\alpha\bx = 0$ say that $[M\bx^\sigma] \in A$ and $[\bx]\in A$ respectively. So, if $A$ contains a point $p$ then it also contains the line of $\mathfrak{S}_\varphi$ through $p$. However $A$, which is a hyperplane of $\PG(V)$, meets every line of $\mathfrak{S}_\varphi$ non-trivially. Hence $A$ contains all lines of $\mathfrak{S}_\varphi$. This canont be, since $\mathfrak{S}_\varphi$ covers $\PG(V)$. So, (\ref{inv3}) holds true, as claimed.  

Let $f^*$ be the semi-linear mapping of $V^*$ which maps every $\xi\in V^*$ onto $\xi^\sigma M^{-1}$ ($\equiv \xi^\sigma M^\sigma$ by (\ref{inv2})) and let $\varphi^* = [f^*]$. Condition (\ref{inv3}) implies that $\varphi^*$ fixes no point of $\PG(V^*)$ while (\ref{inv2}) forces $\varphi^*$ to be an involution. Hence $\varphi^*$ defines a line-spread $\mathfrak{S}^*_{\varphi^*}$ of $\PG(V^*)$. We have $f^*(\xi)\bx = \xi^\sigma M^{-1}\bx = \xi^\sigma M^\sigma\bx = (\xi M\bx^\sigma)^\sigma$ by (\ref{inv2}). Therefore $f^*(\xi)\bx = 0$ if and only if $\xi f(\bx) = 0$. This shows that $\mathfrak{S}^*_{\varphi^*}$ is indeed the dual of $\mathfrak{S}_\varphi$. \hfill $\Box$ 

\begin{note}
\em
Let $\mathfrak{S}_\varphi$ be a semi-standard line-spread and let $\sigma$ be the automorphism of $\FF$ associated to $\varphi$. The hyperplane $\cH_{\mathfrak{S}_\varphi}$ does not arise from $\veS$, since $\mathfrak{S}_\varphi$ is non-standard. However it arises from the following twisted version of $\veS$, introduced by Thas and Van Maldeghem \cite{TVM} for the case $n = 2$ and generalized to any $n \geq 2$ by De Schepper, Schillewaert and Van Maldeghem \cite{DSSVM}:
\[\begin{array}{ccccc}
\ve_\sigma & : & A_{1,\{1,n\}}(\FF) & \rightarrow & \PG(M_{n+1}(\FF)) \\
 & & ([\bx],\xi]) & \rightarrow & [\bx^\sigma\otimes\xi]
\end{array}\]  
It is natural to call the hyperplanes of $A_{1,\{1,n\}}(\FF)$ which arise from embeddings as $\ve_\sigma$ but not from $\veS$ {\em hyperplanes of twisted type}. We are not going to discuss these hyperplanes in this paper. We keep them for future investigations. 
\end{note} 

\section{Proof of Propositions \ref{hyp pure} and \ref{main1}}\label{pure}

{\bf Proof of Proposition \ref{hyp pure}.} For a non-zero vector $\ba\in V$, the points of ${\cal M}_{[\ba]}$ correspond via $\veS$ to points of $\PG(M_{n+1}(\FF))$ represented by pure tensors $\ba\otimes \xi$ with $\xi(\ba) = 0$. Given a matrix $M\in M_{n+1}(\FF)\setminus\langle I\rangle$, we have ${\cal M}_{[\ba]}\subseteq\cH_M$ if and only if $\xi(\ba) = 0$ implies $f(M, \ba\otimes \xi) = 0$ for every $\xi\in V^*$. Namely, with $\ba = (a_i)_{i=0}^n$ and $M = (m_{i,j})_{i,j=0}^n$,   
\[\sum_{i=0}^n\xi_ia_i = 0 ~ \Rightarrow ~ \sum_{i,j=0}^nm_{i,j}a_j\xi_i = 0, \hspace{5 mm} \forall \xi = (\xi_i)_{i=0}^n\in V^*. \]
This implication amounts to say that $M\ba \in \langle\ba\rangle$, namely $\ba$ is a right eigenvector of $M$. So, claims (1) and (3) of Proposition \ref{hyp pure} are equivalent. In a similar way we see that, for $\alpha\in V^*\setminus\{0\}$, the hyperplane $\cH_M$ contains ${\cal M}_{[\alpha]}$ if and only if $\alpha M\in\langle\alpha\rangle$, namely $\alpha$ is a left eigenvector of $M$. Hence (2) and (3) of Proposition \ref{hyp pure} are equivalent.   \hfill $\Box$ \\  

\noindent
{\bf Proof of Proposition \ref{main1}.} Given a point $a$ and a hyperplane $A$ of $\PG(n,\FF)$, possibly $a\not\in A$, let ${\cal C}_{a,A}$ be the set of points of $A_{n,\{1,n\}}(\FF)$ which are collinear with at least one point of ${\cal M}_a\cup{\cal M}_A$. Note that if $a\in A$ then ${\cal C}_{a,A}$ is the set of points at distance at most $2$ from $(a,A)$. Proving Proposition \ref{main1} amounts to prove that ${\cal C}_{a,A} = \cH_{a,A}$. 

Let $(x,X)\in {\cal C}_{a,A}$. To fix ideas, let $(x,X)\perp (a,K)$ for some point $(a,K)\in{\cal M}_a$. Then either $x = a$ or $X = K$. In both cases $a\in X$. Therefore, if $\ba\in V$ represents $a$ and $\xi\in V^*$ represents $X$, we have $\xi(\ba) = 0$. So, if $\bx$ and $\alpha$ represent $x$ and $A$ respecticely, we have $f(\ba\otimes A, \bx\otimes X) = \xi(\ba)\alpha(\bx) = 0$, namely $(x,X)\in \cH_{a,A}$. By a similar argument, if $(x,X)^\perp\cap{\cal M}_A\neq\emptyset$ then $(x,X)\in \cH_{a.A}$. Thus, we have proved that ${\cal C}_{a,A} \subseteq \cH_{a,A}$. 

Conversely, let $(x,X)\in \cH_{a,A}$. Then $\alpha(\bx)\xi(\ba) = 0$ with $\ba$, $\alpha$, $\bx$ and $\xi$ as above. Therefore either $\xi(\ba) = 0$ or $\alpha(\bx) = 0$. If $\xi(\ba) = 0$ then $(x,X)$ is collinear with a point of ${\cal M}_a$ and if $\alpha(\bx) = 0$ then $(x,X)$ is collinear with a point of ${\cal M}_A$. The reverse inclusion $\cH_{a,A}\subseteq {\cal C}_{a,A}$ is proven. So, ${\cal C}_{a,A} = \cH_{a,A}$. \hfill $\Box$ \\

The next corollary is a trivial consequence of the last claim of Proposition \ref{main1}. It highlights the geometrical meaning of the relation $\perp_f$ between null-traced pure tensors.  

\begin{cor}\label{prop1}
For two pure tensors $\ba\otimes\alpha$ and $\bfb\otimes\beta$ of $(V\otimes V^*)_0$, we have $\ba\otimes\alpha\perp_f \bfb\otimes\beta$ if and only if the corresponding points  $([\ba], [\alpha])$ and $([\bfb], [\beta])$ of $A_{n,\{1,n\}}(\FF)$ are at distance at most $2$ in the collinearity graph of $A_{n,\{1,n\}}(\FF)$. 
\end{cor}

\section{Proof of Lemma \ref{main2 lem} and Theorem \ref{main2}}\label{max} 

{\bf Proof of Lemma \ref{main2 lem}.} Throughout this proof we denote by $\cal G$ the collinearity graph of $A_{n,\{1,n\}}(\FF)$ and by $d((x,X),(y,Y))$ the distance in $\cal G$ between two points $(x,X)$ and $(y,Y)$ of $A_{n,\{1,n\}}(\FF)$. 

Given a point $(a,A)$ of $A_{n,\{1,n\}}(\FF)$, let $\cH = \cH_{a,A}$ be the singular subspace with $(a,A)$ as its deepest point and let ${\cal K}$ be the complement of $\cH$ in the point-set of $A_{n,\{1,n\}}(\FF)$. Recall that $\cal K$ is the set of points $(x,X)$ of $A_{n,\{1,n\}}(\FF)$ such that $d((a,A),(x,X)) = 3$; equivalently, $x\not\in A$ and $a\not\in X$. We shall prove that any two points $(b,B)$ and $(c,C)$ of $\cal K$ can be joined by a path of $\cal G$ all nodes of which are in $\cal K$. The maximality of $\cH$ then follows from Shult \cite[Lemma 4.1.1]{S}.

If $(b,B)\perp(c,C)$ there is nothing to prove. Suppose that $(b,B)\not\perp (c,C)$. So, $b\neq c$ and $B\neq C$. Put $L = B\cap C$. Suppose firstly that $A \not\supseteq L$. Then we can choose a point $d\in L\setminus A$. As $d\not \in A$ and $a\not\in B\cup C$ (because $d((b,B),(a,A)) = d((c,C),(a,A)) = 3$), both points $(d,B)$ and $(d,C)$ belong di $\cal K$. Moreover, $(b,B)\perp (d,B)\perp (d,C)\perp (c, C)$. The path $((b,B), (d,B), (d,C), (c,C))$ has the required properties. Similarly, if $a\not \in \ell = \langle b, c\rangle$, then we can choose a hyperplane $D$ of $\PG(n,\FF)$ through $\ell$ which does not contain $a$. Then $(b,D), (c,D) \in {\cal K}$ and $(b,B)\perp(b,D)\perp(c,D)\perp(c,C)$. The path $((b,B), (b,D), (c,D), (c,C))$ does the job.   

With $L$ and $\ell$ as above, suppose now that $a\in \ell$ and $L \subset A$. Clearly, neither $b$ nor $c$ belong to $L$. Choose a point $d\in B\setminus(L\cup\{b\})$. Then $(d,B) \in {\cal K}$ and $(b,B)\perp(d,B)$. Moreover, the line of $\PG(n,\FF)$ through $d$ and $c$ does not contain $a$. By the previous pragraph, we can join $(d,B)$ with $(c,C)$ in $\cal K$ by a path of lenght 3. Pasting this path with the edge $\{(b,B), (d,B)\}$ we get a path of lenght 4 from $(b,B)$ to $(c,C)$ fully contained in $\cal K$. \hfill $\Box$\\

\noindent  
{\bf Proof of Theorem \ref{main2}.} Throughout this proof we denote the points of $A_{n,\{1,n\}}(\FF)$ by low case german letters, thus avoiding the slightly awkward notation $(a,A)$ used so far. 

For a contradiction, let $\cH_1$ be a hyperplane of $A_{n,\{1,n\}}(\FF)$ properly contained in a proper subspace $\cH_2$ of $A_{n,\{1,n\}}(\FF)$. Clearly, $\cH_2$ is also a hyperplane. 

Choose a point $\gA = (a,A)\in\cH_2\setminus \cH_1$ and let $M$ be a line of $A_{n,\{1,n\}}(\FF)$ through $\gA$. Then $M$ meets $\cH_1$ in a point $\gB$, necessarily different from $\gA$. Hence $M\subseteq \cH_2$. Let $N$ be another line of $A_{n,\{1,n\}}(\FF)$ meeting $M$ in a point $\gC \neq \gA, \gB$. The line $N$ also meets $\cH_1$ in a point, necessarily different from $\gC$ (because $\gC\not \in \cH_1$). Hence $N$ too is fully contained in $\cH_2$. So, $\cH_2$ contains all points of the singular hyperplane $\cH_\gA = \cH_{a,A}$ except possibly the points $\gX$ at distance $2$ from $\gA$ but such that  $\gX^{\perp}\cap \gA^{\perp}\subseteq \cH_1$. Let $\cal X$ be the set of these points.  

Let $\gX\in {\cal X}$ and suppose that $\gX$ and $\gA$ form a polar pair. So $\gX^{\perp}\cap\gA^{\perp} = \{\gB,\gB'\}\subset \cH_1$ for a unique non-collinear pair $\{\gB, \gB'\}$. Let $\Gamma$ be the symp containg $\gA$ and $\gX$. Recall that $\Gamma$ is a grid. The grid $\Gamma$ contains the ordinary quadrangle formed by $\gA, \gB, \gB'$ and $\gX$. Moreover, as $\cH_2$ contains all points of $\cH_\gA$ but possibly those of ${\cal X}$, the intersection $\Gamma\cap\cH_2$ contains all points of $\Gamma$ but possibly $\gX$. Indeed all points of $\Gamma$ different from $\gX$ are collinear with either $\gA$ or a point collinear with $\gA$ and different from both $\gB$ and $\gB'$, hence all of them belong to $\cH_2$. However $\cH_2\cap\Gamma$ is a hyperplane of $\Gamma$. Hence it cannot be the same as the singleton $\{\gX\}$. Consequently, $\gX\in \cH_2$. 

Therefore the points of $\cal X$ exterior to $\cH_2$, if they exist, form special pairs with $\gA$. So, let $\gX = (x,X)$ be a point of $\cal X$ forming a special pair with $\gA$. Then either $a\in X$ but $x\not\in A$ or $x\in A$ but $a\not\in X$. Assume that $a\in X$ and $x\not\in A$. Let $\ell = \langle a, x\rangle$ be the line of $\PG(n,\FF)$ through $a$ and $x$ and let $Y$ be a hyperplane of $\PG(n,\FF)$ containing $\ell$. The point $\gA_Y = (a,Y)$ is collinear with $\gA$. Hence $\gA_Y\in\cH_2$. Moreover, $\gA_Y$ has distance at most $2$ from $\gX$ and, if $d(\gA_Y, \gX) = 2$ then $\{\gA_Y, \gX\}$ is a polar pair. Suppose firstly that $\gA_Y \not \in \cH_1$ for some choice of $Y\supseteq \ell$. Then we can replace $\gA$ with $\gA_Y$ and the same argument previously used for $\gA$ yields $\gX\in \cH_2$. Similarly, if $x\in A$ but $a\not\in X$ and $(y,A) \not \in \cH_1$ for some point $y\in A\cap X$, then $\gX \in \cH_2$. Accordingly, the points at distance $2$ from $\gA$ which might not belong to $\cH_2$ are contained in a proper subset of $\cal X$, which we shall describe in a few lines. 

Recall that $\gA^\perp$ is the union ${\cal M}_a\cup{\cal M}_A$ of the maximal singular subspaces of $A_{n,\{1,n\}}(\FF)$ based at $a$ and $A$ and that ${\cal M}_a\cap{\cal M}_A = \{\gA\}$. We know that $\gA^\perp\subseteq\cH_2$. Moreover, since $\cH_1$ is a hyperplane of $A_{n,\{1,n\}}(\FF)$ but it does not contain $\gA$, the intersections $\cH_1\cap{\cal M}_a$ and $\cH_1\cap{\cal M}_A$ are hyperplanes of ${\cal M}_a$ and ${\cal M}_A$ respectively and neither of them contains $\gA$. Accordingly, there exists a line $\ell$ of $\PG(n,\FF)$ containing $a$ but not contained in $A$ and such that ${\cal M}_a\cap\cH_1 = \{(a,X) \mid X\supseteq \ell\}$. Dually, there exists a sub-hyperplane $L$ of $\PG(n,\FF)$ contained in $A$ but not containing $a$ and such that ${\cal M}_A\cap\cH_1 = \{(x,A)\mid x\in L\}$. The points at distance $2$ from $\gA$ which do not belong to $\cH_2$, if they exist, belong to the following  proper subset of $\cal X$:
\[{\cal X}' ~:= ~ \{(x,X)\in {\cal X} \mid \mbox{ either } X\supset \ell \mbox{ but } x\not\in A \mbox{ or } x\in L \mbox{ but } a\not\in X\}.\] 
Let $\gB = (b,B) \in {\cal X}'$. So, either $B\supset\ell$ and $b\not\in A$ or $b\in L$ but $a\not\in B$. To fix ideas, let $B\supset \ell$ and $b\not\in A$. Suppose that $b\not\in\ell$ and let $m$ be the line of $\PG(n,\FF)$ through $a$ and $b$. A sub-hyperplane $M$ of $\PG(n,\FF)$ exists which is contained in $B$, contains $m$ but does not contain $\ell$. Clearly, $B$ is the unique hyperplane of $\PG(n,\FF)$ which contains both $m$ and $\ell$. Therefore, no hyperplane of $\PG(n,\FF)$ containing $M$ and different from $B$ contains $\ell$. Hence all points $\gB_X = (b,X)$ with $M\subset X\neq B$ belong to $\cH_2$. However, the set of points $(b,X)$ with $X\supset M$ is a line of $A_{n,\{1,n\}}(\FF)$ and $\cH_2$, which is a subspace of $A_{n,\{1,n\}}(\FF)$, contains all points of it but possibly $\gB$. Consequently this line is fully contained in $\cH_2$. Therefore $\gB\in \cH_2$. 

Thus we have proved that if $B\supset\ell$ then all points of ${\cal M}_B$ belong to $\cH_2$ but at most the points $(x,B)$ with $a \neq x\in \ell$. So, ${\cal M}_B\cap\cH_2$ contains the complement of a line of the projective space ${\cal M}_B$ and at least one point of that line. Since ${\cal M}_B\cap\cH_2$ is a subspace of ${\cal M}_B$, this forces ${\cal M}_B\cap\cH_2 = {\cal M}_B$, namely ${\cal M}_B\subseteq \cH_2$. Similarly, if $b\in L$ then ${\cal M}_b\subseteq \cH_2$. So, ${\cal X}'\subseteq \cH_2$. It follows that $\cH_2\supseteq \cH_\gA$. However $\cH_\gA$ is a maximal subspace of $A_{n,\{1,n\}}(\FF)$ by Lemma \ref{main2 lem}. Hence $\cH_2 = \cH_\gA$.

The point $\gA$ considered so far is an arbitrary point of $\cH_2\setminus \cH_1$. Therefore $\cH_2 = \cH_\gA$ for every point $\gA\in\cH_2\setminus\cH_1$. Accordingly, $\cH_\gA = \cH_\gB$ for any choice of $\gA, \gB\in \cH_2\setminus\cH_1$. However $\cH_2\setminus\cH_1$ contains more than one point. Hence different points $\gA$ and $\gB$ exist such that $\cH_\gA = \cH_\gB$. This conclusion contradicts the uniqueness of the deepest point of a singular hyperplane.  \hfill $\Box$  

\section{Proof of Theorems \ref{main3} and \ref{main4}}\label{spread} 

{\bf Proof of Theorem \ref{main3}.}  Recall that $A_{n,\{1,n\}}(\FF)$ admits two types of lines, namely the lines $L_{p.L} = \{(p,H) \mid H\supset L\}$ for $L$ a sub-hyperplane of $\PG(n,\FF)$ and $p\in L$ and the lines $L_{\ell, H} = \{(p,H) \mid p\in \ell\}$ for a line $\ell$ and a hyperplane $H$ of $\PG(n,\FF)$ containing $\ell$.  

Consider $L_{p,L}$. If $L\in\mathfrak{S}^*$ then $L_{p,L}\subseteq\cH_\mathfrak{S}$, by definition of $\cH_\mathfrak{S}$ and Lemma \ref{main3 lem3}. Otherwise, let $\ell_p$ be the line of $\PG(n,\FF)$ through $p$ and put $H_0 = \langle \ell_p, L\rangle$. Then $(p,H_0)$, which is a point of $L_{p,L}$, belongs to $\in\cH_\mathfrak{S}$, by the definition of $\cH_\mathfrak{S}$. If $(p,H)$ is another point of $L_{p,L}$, different from $(p,H_0)$, then $H$ cannot contain $\ell_p$. Hence $(p,H)\not\in\cH_\mathfrak{S}$. So, $L_{p,L}\cap\cH_\mathfrak{S} = \{(p,H_0)\}$. 

Similarly, given a line $L_{\ell,H}$, if $\ell\in\mathfrak{S}$ then $L_{\ell, H}\subseteq \cH_\mathfrak{S}$. In contrast, let $\ell\not\in \mathfrak{S}$ and let $L_H$  be the unique member of $\mathfrak{S}$ contained in $H$. Then $\ell$ meets $L_H$ in a unique point, say $p_0$, and $L_{\ell,H}\cap\cH_\mathfrak{S} = \{(p_0,H)\}$.  \hfill $\Box$ \\

\noindent
{\bf Proof of Theorem \ref{main4}.} Suppose that $M = (m_{i,j})_{i,j=0}^n$ admits no eigenvalue in $\FF$ and satisfies conditions $(\mathrm{S}_{\mathrm{mat}})$ and $\mathrm{S}_{\mathrm{mat}}^*)$ of Lemma \ref{main4 lem1} and let $\mathfrak{S}_M$ and $\mathfrak{S}^*_M$ be as in Lemma \ref{main4 lem2}. We shall prove that $\cH_{\mathfrak{S}_M} = \cH_M$. 

Indeed, for a pure tensor $\bx\otimes \xi\in (V\otimes V^*)_0$, we have $\bx\otimes \xi\perp_f M$ if and only if $\sum_{i,j}m_{i,j}x_j\xi_i = 0$, where $(x_i)_{i=0}^n = \bx$ and $(\xi_i)_{i=0}^n = \xi$. Equivalently, $\xi M\bx = 0$, namely $\xi(M\bx) = 0$. On the other hand, $\xi(\bx) = 0$ because $\bx\otimes \xi\in (V\otimes V^*)_0$. Therefore the hyperplane $[\xi]$ of $\PG(n,\FF)$ represented by $\xi$ contains both $[\bx]$ and $[M\bx]$, namely $[\xi]$ contains the line $\langle [\bx], [M\bx]\rangle$, which is indeed the line of $\mathfrak{S}_M$ through $[\bx]$. The equality $\cH_M = \cH_{\mathfrak{S}_M}$ is proven.  

Conversely, given $M \in M_{n+1}(\FF)\setminus\langle I\rangle$, suppose that $\cH_M$ is of spread type. A hyperplane of spread type contains no maximal singular subspace of $A_{n,\{1,n\}}(\FF)$. Hence $M$ admits no eigenvalue in $\FF$, by  Proposition \ref{hyp pure}. Moreover, $\cH_M = \cH_{\mathfrak{S}_N}$ for a matrix $N\in M_{n+1}(\FF)$ satifying property $(\mathrm{S}_{\mathrm{mat}})$ and admitting no eigenvalues in $\FF$. However $\cH_{\mathfrak{S}_N} = \cH_N$ by the previous paragraph. Hence $\cH_M = \cH_N$. Therefore $M = sN + tI$ for suitable scalars $s, t\in \FF$, $s\neq 0$, by Remark \ref{rem notation}. We have $N^2\bx \in \langle \bx, N\bx\rangle$ for every $\bx\in V$, by property $(\mathrm{S}_{\mathrm{mat}})$ on $N$. So, $M^2\bx = s^2N^2\bx + 2st N\bx + t^2\bx \in \langle \bx, N\bx\rangle$. Moreover $\langle \bx , N\bx\rangle = \langle \bx, M\bx\rangle$, since $N = s^{-1}M -  s^{-1}tI$. It follows that $M^2\bx\in \langle \bx, M\bx\rangle$. So, $M$ also satisfies  $(\mathrm{S}_{\mathrm{mat}})$ (equivalently, $(\mathrm{S}_{\mathrm{mat}}^*)$).  \hfill $\Box$

\end{document}